# TRANSFORMATIONS, DYNAMICS AND COMPLEXITY

## Nikolaj M. Glazunov

*National Aviation University,
03680, Kiev-680 GSP Ukraine,
glanm@yahoo.com*

Abstract: We review and investigate some new problems and results in the field of dynamical systems generated by iteration of maps, $\beta$-transformations, partitions, group actions, bundle dynamical systems, Hasse-Kloosterman maps, and some aspects of complexity of the systems.

## INTRODUCTION

Dynamical systems have several important complexity measures among which are measure-theoretic, topological and metric entropies. A. Kolmogorov defines the complexity of a measure-preserving transformation by generators. A generator for a measure-preserving transformation $T$ is a partition $\xi$ with finite entropy such that the set of finite entropy partitions subordinate to some $V_{i=-n}^{n} T^{-i}(\xi)$ is dense in the set of finite entropy partitions endowed with the Rokhlin metric. For smooth dynamical systems topological entropy characterizes the total exponential complexity of the orbit structure. Metric entropy with respect to an invariant measure codes the exponential growth rate of the statistically significant orbits. In symbolic dynamics complexity of a dynamical system is measured with respect to a coding of its orbits. We review and investigate some new problems and results in the field of dynamical systems generated by iteration of maps, $\beta$-transformations, partitions, group actions, bundle dynamical systems, Hasse-Kloosterman maps, and some aspects of complexity of the systems.

## 2. SEQUENCES, THEIR DYNAMICS AND COMPLEXITY

*2.1 Iteration of the "$\times 2$ maps" $T(x) = 2x \pmod{1}$.*

Let $\mathcal{M}$ be the simplex of × 2-invariant Borel probability measures on $X = [0,1]$, $b(\mu) = \int x d\mu(x)$ the barycentre of the probability measure $\mu$, $\mathcal{M}_\rho = \{\mu \in \mathcal{M}: b(\mu) = \rho\}$. In paper (Jenkinson, 2009) Jenkinson proved

**Theorem** (Jenkinson, 2009). For every $\varrho \in [0,1]$, the ordered set $(\mathcal{M}_\varrho, \prec)$ has a least element. This least element is the Sturmian measure $S_\varrho$ of rotation number $\varrho$.

The theorem itself is a corollary of a result on $C^2$ convex functions: let $f: X \to \mathbf{R}$ be a $C^2$ convex function. For every $\varrho \in (0,1)$ there exists $\theta \in \mathbf{R}$ such that the Sturmian measure $S_\varrho$ is a minimizing measure for the function $f_\theta$ defined by the

$$f_\theta(x) = f(x) + \theta x.$$

A number of corollaries about Sturmian measures in $\mathcal{M}_\varrho$ are also given, including facts on smaller $f$ −integral, on smallest variance and the result that Sturmian orbit has largest geometric mean. The result on $C^2$ convex functions comes partially from a reformulation of precondition de Sturm by T. Bousch (2000). The paper (Jenkinson, 2008) includes a short survey of ergodic optimization and of Sturmian measures as minimizing and maximizing measures for classes of functions. If $\mathcal{M}$ is the set of Borel probabilities on the unit circle which are invariant with respect to the mapping $T(x) = 2x \ (mod \ 1)$ then (Bousch, 2000) proved that for each function $\rho_\omega(t) \coloneqq \cos(2\pi(t-\omega))$, there is exactly one element $\mu \in \mathcal{M}$ which maximizes $\int \rho_\omega d\mu$, and that the support of this measure is contained in a semicircle. In particular, the image of the map $\mu \mapsto \int \exp\{2i\pi t\} d\mu$ is compact and a convex set of the complex plane which does not contain any line segments on its boundary. T. Bousch (2000) also proved that the maximizing measure is periodic for every $\omega$ except on a set which has measure zero and Hausdorff dimension zero.

*2.2 Baker's transformations, β-transformations and their extensions.*

For $\beta > 1$ and $x \in [0, 1]$ let

$$T_\beta(x) = \beta x \ mod \ 1,$$

be the β-transformation and

$$d_\beta(x) = (x_i)_{i \geq 1}, x_i = [\beta T_\beta^{i-1}(x)],$$

the β-expansion. The symbolic dynamics of the β-transformation has investigated by A. Rénei who introduced the transformation and proved that it is ergodic (Rénei, 1957) and by W. Perry who described possible sequences that can be a β-expantion (Perry, 1960). F. Blanchard is noticed the key role of $d_\beta(1)$ (Blanchard, 1989). J. Alexander and J. Yorke and C. Bose have investigated respectively fat baker's transformations (Alexander and Yorke, 1984) and generalized baker's transformations (Bose, 1989). The natural extensions of β-transformations have been characterized in paper ( Dajani *et al.*, 1996 ). In paper (Broun and Yin, 2000) the authors study the set $\mathcal{S}$ of all those β for which the natural extension of $T_\beta$ can be represented by the map

$$\mathcal{T}_\beta(x,y) = \left(T_\beta x, \frac{[\beta x] + y}{\beta}\right)$$

defined on a simply connected subset of the unit square, and with invariant measure a constant multiple of the 2-dimensional Lebesgue measure on $[0,1]^2$. In paper (Kwon, 2009) the author investigates sufficient conditions for the existence of β-transformations such that the natural extensions of the transformations can be viewed as generalized baker's transformations. He characterizes such $\beta$ and studies their properties. This paper includes recent results on the symbolic dynamics of the transformations.

## 3. RETURN AND HITTING TIMES IN DYNAMICAL SYSTEMS GENERATED BY A TYPICAL PARTITIONS

In the paper (Grzegorek and Kupsa, 2009) authors investigate limit distributions for the return and hitting times appearing in $\alpha$-mixing processes and substantial limit distributions for the return times in adding machines. They state new results on limit distributions.

It was shown by T. Downarowicz and Y. Lacroix and by T. Downarowicz, Y. Lacroix and D. Leandri that for every non-periodic ergodic system, the zero function is a limit distribution of the hitting times to cylinders achieved along a subsequence of lengths of upper density 1 and with probabilities increasing to 1, in every typical process defined on this system. In the present paper it is shown that the exponential distribution

$$E(t) = max(0, 1 - exp(-t))$$

is another limit distribution (for both the return and hitting times) achieved in Downarowicz and Lacroix manner in every typical process generated on any $\alpha$-mixing system. Authors of the paper (Grzegorek and Kupsa, 2009) generalize results obtained by M. Abadi and by A. Galves and B. Schmitt for $\phi$-mixing and exponentially $\alpha$-mixing processes and prove in their theorem 4 that for an $\alpha$-mixing process the exponential distribution is the unique limit distribution for the hitting time. In the last section of (Grzegorek and Kupsa, 2009) authors show that for an adding machine the piecewise linear function is a substantial limit distribution for the return times in a process generated by a typical partition (theorem 2).

## 4. MEASURE-THEORETIC ENTROPY AND BUNDLE DYNAMICAL SYSTEMS

A. Kolmogorov defines the complexity of a measure-preserving transformation by generators. A generator for a measure-preserving transformation $T$ is a partition $\xi$ with finite entropy such that the set of finite entropy partitions subordinate to some $\bigvee_{i=-n}^{n} T^{-i}(\xi)$ is dense in the set of finite entropy partitions endowed with the Rokhlin metric. Random dynamical systems relate a partial case of bundle dynamical systems by I. Cornfeld, S. Fomin, and Ya. Sinaì (Корнфельд *et al.,* 1980). In paper (Zhu, 2009) the author investigates measure-theoretic entropy of random dynamical systems (RDS) and presents in the framework of RDS the extension of local entropy formula by M. Brin and A. Katok and the definition of the measure-theoretic entropy using spanning set by A. Katok.

The main results of the paper (Zhu, 2009) are local entropy formula and the presentation of the measure-theoretic entropy of random dynamical systems using spanning set. The last theorem is an extension of the result presented by M. Pollicott.

The proof of the local entropy formula uses Shannon-McMillan-Breiman Theorem by T. Bogenschutz and Egorov Theorem by P. Halmos. In the proof of last theorem the author (Kwon, 2009) is based on the proof of the variational principle by M. Misiurewicz.

## 5. AMENABILITY, NON AMENABILITY AND COMPLEXITY

A countable discrete group $G$ is called amenable if it has an invariant mean or equivalently if it has a left invariant finitely additive measure $\mu, \mu(G) = 1$. Amenability of groups is introduced by J. von Neumann (1929) (he called such groups "measurable") to explain results of the paper by S. Banach and A. Tarski (1924) (the so-called Banach-Tarski paradox).

There are at least twelve equivalent definitions of amenability. Von Neumann proved that a group is not amenable provided it contains a free subgroup on two generators. First counterexamples to the von Neumann problem were constructed by A. Ol'shanskii (1980) (the so called "Tarski monsters"). He proved that the groups with all proper subgroups cyclic constructed by him, both torsion-free and torsion, are not amenable. A finitely presented non-amenable group without free non-cyclic subgroups have constructed by A. Ol'shanskii and M. Sapir (2003).

## 6. FUNCTIONAL CASE OF THE HASSE-KLOOSTERMAN MAP

Let $F_p$ be the prime finite field, $F_p^*$ multiplicative group of $F_p$. Below in the section we follow to (Glazunov, 2010).

### 6.1. *Hasse component.*

Let
$$y^2 = f(x), \ f(x) = x^3 + cx + d,$$

be a cubic polynomial in prime finite field $F_p$. For the number $\#C_p$ of points of the curve $C: y^2 = f(x)$ in $F_p$ the well known formula

$$\#C_p = \sum_{x=0}^{p-1}\left(1+\left(\frac{f(x)}{p}\right)\right)$$

take place, where $\left(\frac{f(x_0)}{p}\right)$ is the Legendre symbol with a numerator which is equal to the value of the polynomial $f(x_0)$ in point $x_0 \in F_p$. It is ease to see that $\#C_p = p - a_p$, where

$$a_p = -\sum_{x=0}^{p-1}\left(\frac{f(x)}{p}\right).$$

If $C$ is the elliptic curve $E$, then the number of points $\#E(F_p) = \#E_p$ of the projective model of the curve $E$ in $F_p$ is represented by the formula $\#E_p = 1 + p - a_p$, where $a_p = 2\sqrt{p} \cdot \cos\varphi_p$. If $C$ is not the elliptic curve, then the value $a_p$ is equal 1, -1 or 0 and ease to compute. In both cases compute: $\varphi_p = \arccos\left(\frac{a_p}{2\sqrt{p}}\right)$ and reduce it to the interval $[0, \pi]$.

### 6.2. *Kloosterman component.*

Let $cd \not\equiv 0 \mod p$,

$$T_p(c,d) = \sum_{x=1}^{p-1} e^{2\pi i \frac{\left(cx+\frac{d}{x}\right)}{p}},$$

the Kloosterman sum. By A. Weil,

$$T_p(c,d) = 2\sqrt{p}\cos\theta_p(c,d).$$

Compute $T_p$, $\cos\theta_p$, $\theta_p$ and reduce $\theta_p$ to the interval $[0, \pi]$.

### 6.3. *Hasse-Kloosterman map.*

Functional case of the Hasse-Kloosterman (HK) map is defined on $(F_p)^* \times (F_p)^*$ with values in $\Pi = [0, \pi] \times [0, \pi]$ and has the form

$$hk(c,d) = (\varphi_p(c, d), \theta_p(c, d)).$$

As *c,d* runs independently the multiplicative group $(F_p)^*$, their product is not divided by *p*, so the map *hk(c,d)* is defined in all points.

*6.4. Coding.*

Let $\mathcal{R}_1$ and $\mathcal{R}_2$ be two finite partitions of the same cardinality *d* of the interval $[0, \pi]$. We call $\mathcal{R}_1$ horizontal, $\mathcal{R}_2$ vertical partition of $[0, \pi]$ and the pair $(\mathcal{R}_1, \mathcal{R}_2)$ *p*-pair. Denote elements of a partition by integer numbers *0,1, ..., d-1*. Functional mapping HK is coded by finite sequence $b_s b_{s-1}...b_1 b_0.a_1 a_2...a_r$. A value of the sequence may be interpreted as a rational number *(x, y)* from unit square, if we put $x = \sum_{i=1}^{r} a_i/d^i$ (respectively, $= \sum_{i=1}^{s} b_{i-1}/d^i$) for *d*-adic expansion of *x* (respectively, for *y*).

# 7. ARITHMETIC CASE OF THE HASSE-KLOOSTERMAN MAP

Let *C* is an cubic algebraic curve over ring of integers $\mathbb{Z}$, *Spec* $\mathbb{Z}$ the affine scheme over $\mathbb{Z}$, $d_1$ the set of points of *Spec* $\mathbb{Z}$, (which correspond prime *p*) where $c, d \equiv 0 \pmod{p}$, $Sp = Spec\ \mathbb{Z} \setminus d_1$.

*7.1 Hasse component.*

If *C* is an elliptic curve *E* over $\mathbb{Z}$, then *E* has good reduction in all points with the exception of prime divisors of it's discriminant. The number of points $\#E(F_p) = \#E_p$ of the curve *E* under it localization *mod p* is expressed by the formula $\#E_p = 1 + p - a_p$, where $a_p = 2\sqrt{p} \cdot \cos\varphi_p$, and the curve *E* is considered as projective. In the case when the localization of *C mod p* is not the elliptic curve the value of $a_p$ is equal 1, -1 or 0 and ease to compute. In both cases compute: $\varphi_p = \arccos\left(\dfrac{a_p}{2\sqrt{p}}\right)$ and reduce it to the interval $[0, \pi]$.

*7.2 Kloosterman component.*

Let $cd \not\equiv 0 \mod p$,

$$T_p(c,d) = \sum_{x=1}^{p-1} e^{2\pi i \frac{\left(cx + \frac{d}{x}\right)}{p}},$$

the Kloosterman sum. By A. Weil,

$$T_p(c,d) = 2\sqrt{p}\ \cos\theta_p(c,d).$$

Compute $T_p$, $\cos\theta_p$, $\theta_p$ and reduce $\theta_p$ to the interval $[0, \pi]$.

*7.3 Hasse-Kloosterman map.*

Arithmetic mapping of the Hasse-Kloosterman (HK) is defined on the product $Sp \times Sp$ with value in $\Pi = [0, \pi] \times [0, \pi]$ and has the form

$$hk(c,d) = (\varphi_p(c, d), \theta_p(c, d)).$$

As *Sp* by its definition does not contains prime divisors of $cd$, the mapping $hk(c,d)$ is defined in all points.

*7.4. Coding.*

Let $\mathfrak{R}_1$ and $\mathfrak{R}_2$ be two finite partitions of the same cardinality *d* of the interval $[0, \pi]$. We call $\mathfrak{R}_1$ horizontal, $\mathfrak{R}_2$ vertical partition of $[0, \pi]$ and the pair $(\mathfrak{R}_1, \mathfrak{R}_2)$ *p*-pair. Denote elements of a partition by integer numbers $0, 1, \ldots, d-1$. Arithmetic mapping HK is coded by infinite sequence $\ldots b_s b_{s-1} \ldots b_1 b_0 . a_1 a_2 \ldots a_r \ldots$. A value of the sequence may be interpreted as a real number $(x, y)$ from unit square, if we put $x = \sum_{i=1}^{\infty} a_i/d^i$ (respectively, $y = \sum_{i=1}^{\infty} b_{i-1}/d^i$ ) for *d*-adic expansion of *x* (respectively, for *y*).

## 8. CONCLUSIONS

We gave a review and present some new results in the field of dynamical systems generated by iteration of maps, *β*-transformations, Hasse-Kloosterman maps, partitions, bundle dynamical systems and some aspects of complexity of the systems.